\documentclass{article}

\usepackage{arxiv}

\usepackage{amsmath,amssymb,enumerate}
\usepackage{amsthm}
\usepackage[utf8]{inputenc} 
\usepackage[T1]{fontenc}    
\usepackage{hyperref}       
\usepackage{url}            
\usepackage{booktabs}       
\usepackage{amsfonts}       
\usepackage{nicefrac}       
\usepackage{microtype}      
\usepackage{cleveref}       
\usepackage{lipsum}         
\usepackage{graphicx}
\usepackage{natbib}
\usepackage{doi}

\usepackage{multirow}
\usepackage{threeparttable}

\newtheorem{proposition}{Proposition}

\title{Enhancing Quadratic Programming Solvers via Quadratic Nonconvex Reformulation}

\usepackage{authblk}

\newbox{\orcid}
\author[1]{Cheng Lu}

\newbox{\orcid}
\author[1]{Yu Fei}

\newbox{\orcid}
\author[1]{Gaojian Kang}

\newbox{\orcid}
\author[1]{Guangtai Qu}

\newbox{\orcid}
\author[2]{Zhibin Deng}

\newbox{\orcid}
\author[3]{Qingwei Jin}

\newbox{\orcid}
\author[4]{Shu-Cherng Fang}

\affil[1]{School of Economics and Management, North China Electric Power University, Beijing, 102206, China}
\affil[2]{School of Economics and Management, University of Chinese Academy of Sciences, Beijing 100190, China}
\affil[3]{School of Management, Zhejiang University, Hangzhou, 310058, China}
\affil[4]{Department of Industrial and Systems Engineering, North Carolina State University, Raleigh,
	NC 27695-7906, USA}

\renewcommand{\undertitle}{}

\begin{document}
\maketitle

\begin{abstract}
	In this paper, we consider solving nonconvex quadratic programming problems using modern solvers such as Gurobi and SCIP. It is well-known that the classical techniques of quadratic convex reformulation can improve the computational efficiency of global solvers for mixed-integer quadratic optimization problems. In contrast, the use of quadratic nonconvex reformulation (QNR) has not been previously explored. This paper introduces a QNR framework--an unconventional yet highly effective approach for improving the performance of state-of-the-art quadratic programming solvers such as Gurobi and SCIP. Our computational experiments on diverse nonconvex quadratic programming problem instances demonstrate that QNR can substantially accelerate both Gurobi and SCIP. Notably, with QNR, Gurobi achieves state-of-the-art performance on several benchmark and randomly generated instances.
\end{abstract}


\section{Introduction}
We study the following general quadratically constrained quadratic programming problem (QCQP):
\begin{equation}\label{QCQP}
	(\textrm{QCQP})~\left\{ \begin{aligned}
		\min ~&~\frac{1}{2}x^{T}Q_0 x+c_0^{T}x \\
		\mbox{s.t.} ~&~\frac{1}{2}x^{T}Q_i x+c_i^{T}x\leq b_i,~i=1,2,\ldots,m,\\
		~&~ 0\leq x_i\leq 1,~i=1,\ldots,n.
	\end{aligned}
	\right.
\end{equation}
where $Q_i \in \mathbb{R}^{n\times n}$, $c_i\in \mathbb{R}^{n}$  for $i=0,1,\ldots,m$, and $b_i\in\mathbb{R}$  for $i=1,\ldots,m$.

The QCQP is among the most fundamental problems in nonlinear optimization and is NP-hard in general, even when $Q_0$ has only a single negative eigenvalue \citep{Pardalos1991}. A variety of branch-and-bound-based global algorithms have been proposed for different classes of nonconvex quadratic optimization problems \citep{Burer2008,Cambini2005,Elloumi2019,Linderoth,Locatelli2024,Tawarmalani,Vandenbussche2005}. In parallel, several general-purpose global solvers--such as CPLEX, Gurobi, BARON \citep{Sahinidis,Tawarmalani}, SCIP \citep{Bestuzheva2025,Vigerske}, $\alpha$-BB \citep{Adjiman1998b,Adjiman1998}, and GloMIQO \citep{Misener}--have been developed for solving nonconvex quadratic programming problems.

Despite significant advances, these solvers can still struggle with certain classic QCQP instances, even of modest size (e.g., with only tens of variables). To address this limitation, we introduce a quadratic nonconvex reformulation (QNR) framework--an unconventional approach that transforms one nonconvex QCQP into another nonconvex one with a modified structure, thereby tightening the McCormick relaxation \citep{McCormick} that underlies many solvers' lower bounds.

A dominant reformulation strategy in the literature of mixed-integer quadratic optimization is the well-known quadratic convex reformulation (QCR) technique, which perturbs the nonconvex quadratic objective function and constraints into convex ones to strengthen the continuous relaxation of a mixed-integer quadratic optimization problem. Originally developed for binary quadratic optimization \citep{Billionnet2007}, QCR has been extended to general quadratic integer optimization \citep{Billionnet2012, Billionnet2013, Billionnet2017}, QCQPs with integer variables \citep{Billionnet2016}, 0-1 polynomial optimization  \citep{Elloumi}, and quadratic programs with cardinality or semi-continuous constraints \citep{Zheng2014, Wu, Zheng2020}. While highly effective for mixed-integer problems, QCR is not applicable to purely continuous nonconvex QCQPs in the framework of modern solvers such as Gurobi and SCIP, since the lower bounds of continuous nonconvex QCQPs are not determined by their continuous relaxations.

For continuous QCQPs, modern solvers typically employ McCormick relaxations \citep{McCormick}. Our approach builds upon this foundation: instead of solving the problem directly, QNR reformulates it to an equivalent formulation which yields a stronger McCormick relaxation, improving the branch-and-bound process without changing the global optimum. This idea is related to matrix cone decomposition and best-DC decomposition methods  \citep{Zheng2011, Zheng2011b}, which also aim to produce high-quality relaxations. However, the goal of QNR is fundamentally different--it seeks a reformulation that directly accelerates a solver rather than designing a stand-alone relaxation.

To the best of our knowledge, this is the first work to explore nonconvex reformulation as a systematic method for improving the efficiency of a QCQP solver. We demonstrate the effectiveness of QNR through extensive computational experiments, showing substantial performance gains for both Gurobi and SCIP across benchmark and randomly generated test sets.

The remainder of this paper is organized as follows. Section \ref{sec2} reviews the McCormick relaxation and introduces the QNR framework for QCQPs. Section \ref{sec3} presents QNR for linearly constrained quadratic programs. Section \ref{sec4} discusses applications to various nonconvex quadratic optimization problems. Section \ref{sec5} reports computational results, and Section \ref{sec6} concludes this paper.

The following notations are used in this paper: We use $\mathbb{S}^n$ to denote the set of $n\times n$ real symmetric matrices, $\mathbb{X}^n$ to denote the set of $n\times n$ symmetric matrices with zero diagonal entries, and $\mathbb{X}_+^n$ to denote the set of matrices in $\mathbb{X}^n$ with nonnegative entries. $\mathbb{R}^n$ denotes the set of $n$-dimensional real vectors. $\mathbb{R}^n_+$ denotes the set of nonnegative vectors in $\mathbb{R}^n$. For a given vector $c\in \mathbb{R}^n$, we use $\textrm{diag}(c)$ to denote the $n\times n$ diagonal matrix with diagonal elements being the respective entries of $c$. For any two matrices $A$, $B\in\mathbb{S}^n$, we define $A\cdot B := \textrm{trace}(AB)$. We denote by $X\succeq 0$ a matrix $X$ that is positive semidefinite, and $X\succ 0$ a matrix that is positive definite.

\section{Quadratic nonconvex reformulation technique}\label{sec2}

This section presents the general framework of the quadratic nonconvex reformulation (QNR) technique for problem (QCQP). Our goal is to improve the computational efficiency, particularly execution time, of state-of-the-art solvers such as Gurobi and SCIP by tightening their lower bounds. We first review the lower bounding approach these solvers use, and then introduce the QNR technique to strengthen those bounds.

\subsection{The McCormick relaxation}

Modern solvers like Gurobi and SCIP typically solve nonconvex quadratic programs using spatial branch-and-cut algorithms. To obtain tractable relaxations, these solvers typically employ McCormick inequalities \citep{McCormick} to approximate nonconvex quadratic functions with convex ones. Understanding this relaxation is essential for developing reformulations that improve the efficiency of these solvers.

Specifically, by introducing the matrix $X:=xx^T$, the quadratic function  $x^T Q_i x/2+c_i^{T}x$ can be linearized as $Q_i \cdot X/2+c_i^{T}x$. The nonconvex constraint $X:=xx^T$ is then relaxed via McCormick inequalities applied element-wise for each pair of $i,j\in\{1,\ldots,n\}$, $i\neq j$:
\begin{equation}	
	(X_{ij},x_i,x_j)\in\mathcal{M}_{ij}:=\left\{ (X_{ij},x_i,x_j)\,\left|\,
	\begin{aligned}
		&X_{ij}- x_i\leq 0,~X_{ij}- x_j\leq 0\\
		&X_{ij}\geq 0,~X_{ij}-x_i -x_j+1\geq 0
	\end{aligned}
	\right.\right\},
\end{equation}
along with $X_{ii}\geq x_i^2$ and $X_{ii} \leq x_i$ for $i\in\{1,\ldots,n\}$.

Both solvers Gurobi and SCIP begin by detecting convexity in the objective function and constraints. McCormick relaxations are applied only to nonconvex quadratic terms--that is, terms where $Q_i$ is not positive semidefinite. Denote by $\mathcal{P}\subseteq\{1,\ldots,m\}$ the set of indices $i$ such that $Q_i\succeq 0$. The relaxation of (QCQP) can then be expressed as:
\begin{equation}\label{MCR}
	(\textrm{MCR})~\left\{\begin{aligned}
		\min ~&~ f(x) \\
		\mbox{s.t.} ~&~ \frac{1}{2}x^{T}Q_i x+c_i^{T}x\leq b_i,~i\in\mathcal{P},\\
		~&~ \frac{1}{2}Q_i \cdot X+c_i^{T}x\leq b_i,~i\in\{1,\ldots,m\}\setminus \mathcal{P},\\
		~&~ x_i^2\leq X_{ii} \leq x_i,~i=1,\ldots,n,\\
		~&~ (X_{ij},x_i,x_j)\in \mathcal{M}_{ij}, ~i\neq j, ~i,j=1,\ldots,n,
	\end{aligned} \right.
\end{equation}
where $f(x)=x^{T}Q_0 x/2+c_0^{T}x$ if $Q_0\succeq 0$, and $f(x)=Q_0\cdot X/2+c_0^{T}x$ otherwise. Note that under the constraints $X_{ii}\geq x_i^2$ and $X_{ii} \leq x_i$, the bounds $0\leq x_i\leq 1$ become redundant and can be dropped.
This formulation is referred to as the McCormick relaxation of (QCQP). Quadratic equality constraints are also relaxed in this framework when present.

\textbf{Remark 1.} Both Gurobi and SCIP implement multiple lower bounding strategies. The McCormick relaxation is one selectable option. For more details on SCIP's lower bounding methods, see \citep{Bestuzheva2025}. For Gurobi's McCormick relaxation implemented, we may refer to the official document in the website\footnote{\url{https://www.gurobi.com/events/non-convex-quadratic-optimization/}}. Although Gurobi is not open-source, its lower bounding methods for nonconvex QCQP are user-accessible.

\textbf{Remark 2.} We choose Gurobi and SCIP as representative solvers due to their widespread use--commercial and open-source, respectively. In addition to the two solvers, the QNR technique applies broadly to any solver employing McCormick relaxations.

\subsection{A general framework of QNR for (QCQP)}\label{sec22}
Building on the McCormick relaxation, QNR aims to tighten the relaxation bound of (QCQP). The core idea is to decompose each matrix $Q_i$ as $Q_i=(Q_i-Z_i)+Z_i$ $(i=0,\ldots,m)$, where $Z_i\in\mathbb{S}^{n}$ is chosen such that $Q_i-Z_i\succeq 0$. The problem can then be reformulated as follows:
\begin{equation}\label{QNR}
	(\textrm{QNR})~\left\{\begin{aligned}
		\min ~&~\frac{1}{2}x^{T}(Q_0-Z_0)x+ c_0^{T}x + \frac{1}{2}t_0 \\
		\mbox{s.t.} ~&~\frac{1}{2}x^{T}(Q_i-Z_i)x+c_i^{T}x+\frac{1}{2}t_i\leq b_i,~i=1,\ldots,m,\\
		~&~ 0\leq x_i \leq 1,~i=1,\ldots,n,\\
		~&~ t_i= x^{T}Z_i x,~i=0,\ldots,m.
	\end{aligned}\right.
\end{equation}
Crucially to note that all nonconvexity is now isolated in the equality constraints $t_i = x^T Z_i x$, while the objective function and inequality constraints become convex. This structural shift can significantly tighten the McCormick relaxation when the parameters $(Z_0,\ldots,Z_m)$ are chosen appropriately. Consequently, solving (QNR) with modern solvers often yields stronger bounds and improved performance compared to solving the original (QCQP).

This approach relates conceptually to matrix cone decomposition and polyhedral approximations \citep{Zheng2011}, where quadratic functions are similarly decomposed and relaxed. However, to our knowledge, no prior work has exploited such nonconvex reformulations explicitly to accelerate a solver. Leveraging recent advances in solvers such as Gurobi and SCIP, this paper is the first systematic study of this strategy.

\subsection{Best parameters for QNR}

We now discuss how to select the parameters $(Z_0, \ldots, Z_m)$ for (QNR). Let $\theta(Z_0,\ldots,Z_m)$ denote the McCormick relaxation bound of (QNR), defined as the optimal value of the following problem:
\begin{equation}\label{QNR-MCR}
	\left\{ \begin{aligned}
		\min ~&~ \frac{1}{2}x^{T}(Q_0-Z_0)x+c_0^{T}x + \frac{1}{2} Z_0\cdot X \\
		\mbox{s.t.} ~&~ X_{ii}\geq x_i^2,~X_{ii}- x_i\leq 0,~i=1,\ldots,n,\\
		~&~ \frac{1}{2}x^{T}(Q_i-Z_i)x+c_i^{T}x+ \frac{1}{2} Z_i\cdot X\leq b_i,~i=1,\ldots,m,\\	~&~ \left(X_{ij},x_i,x_j\right)\in \mathcal{M}_{ij},~i,j=1,\ldots,n,~i\neq j.
	\end{aligned}\right.
\end{equation}
The parameters are determined by solving the following problem:
\begin{equation}\label{best_par}
	\left\{\begin{aligned}
		\max ~&~ \theta(Z_0,\ldots,Z_m) \\
		\mbox{s.t.:} ~&~ Q_i-Z_i\succeq 0,~i=0,\ldots,m.
	\end{aligned}\right.
\end{equation}
By results in \citep[Theorem 1]{Billionnet2016} and \citep[Theorem 1]{Zheng2011b}, under the assumption that (QCQP) is strictly feasible, the optimal solution $(Z^\ast_0, \ldots, Z^\ast_m)$ of \eqref{best_par} can be recovered as the solution to a semidefinite programming relaxation of (QCQP), summarized below:
\begin{proposition}\label{thm1}
	Assuming (QCQP) is strictly feasible. The optimal solution $(Z^\ast_0,\ldots,Z^\ast_m)$ of \eqref{best_par} is given by
	\begin{equation}\label{recoverZ}
		\left\{
		\begin{aligned}
			&Z^\ast_i=Q_i,~i=1,\ldots,m,\\
			&Z^\ast_0=-\sum_{i=1}^m \gamma^\ast_i Q_i-2\textrm{diag}(\beta^\ast)+2M^\ast+2N^\ast-2R^\ast-2S^\ast,
		\end{aligned}
		\right.
	\end{equation}
	where $(\gamma^\ast,\beta^\ast,M^\ast,N^\ast,R^\ast,S^\ast)$ solves the following semidefinite program (QNR-SDP):
	\begin{equation}\label{QNR-SDP}
		(\textrm{QNR-SDP})~\left\{\begin{aligned}
			\max ~&~ \tau \\
			\textrm{s.t.} ~&~ \begin{bmatrix}
				2\hat{b}-2\tau &\hat{c}^T\\
				\hat{c} ~&\hat{Q} \\
			\end{bmatrix}\succeq 0,\\
			~&~\hat{Q}=Q_0+\sum_{i=1}^m \gamma_i Q_i +2\textrm{diag}(\beta)-2M-2N+2R+2S,\\
			~&~\hat{c}=c_0+\sum_{i=1}^m \gamma_i c_i -\beta+2Ne-Re-Se,\\
			~&~\hat{b}=-\sum_{i=1}^m \gamma_i b_i-\sum_{i=1}^n\sum_{j=1}^n N_{ij},\\
			~&~\beta\in \mathbb{R}^n_+,~\gamma\in \mathbb{R}^m_+, ~\tau\in\mathbb{R}, ~M,N,R,S\in\mathbb{X}_+^n.
		\end{aligned}\right.
	\end{equation}
\end{proposition}

A full proof for Proposition~\ref{thm1} is provided in Appendix A.

After solving (QNR-SDP), the problem (QNR) can be constructed and passed to the solver. We call the step of solving (QNR-SDP) and constructing (QNR) as QNR preprocessing step. This preprocessing step tightens the McCormick relaxation and thus improves the performance of a solver.

Notice that the dual of (QNR-SDP) corresponds to the well-known SDP+RLT relaxation \citep{Anstreicher}:
\begin{equation}\label{SDPRLT}
	(\textrm{SDP+RLT})~\left\{\begin{aligned}
		\min ~&~ \frac{1}{2}Q_0\cdot X+c_0^{T}x \\
		\mbox{s.t.} ~&~ \frac{1}{2}Q_i \cdot X+c_i^{T}x\leq b_i,~i=1,\ldots,m,\\
		~&~ X_{ii}- x_i\leq 0,~i=1,\ldots,n,\\
		~&~ (X_{ij},x_i,x_j)\in \mathcal{M}_{ij},~i,j=1,\ldots,n,~i\neq j,\\
		~&~ X\succeq xx^T.
	\end{aligned}\right.
\end{equation}
Thus, under strict feasibility of (QCQP), the root-node bound of (QNR) matches the SDP+RLT bound when using the parameters from Proposition~\ref{thm1}.

\section{QNR for linearly constrained quadratic programming problems}\label{sec3}
In this section, we focus on the linearly constrained quadratic programming problem defined as follows:
\begin{equation}\label{LCQP}
	(\textrm{LCQP})~\left\{\begin{aligned}
		\min ~&~\frac{1}{2}x^{T}Q x +c^Tx\\
		\mbox{s.t.} ~&~ a_i^T x = d_i,~i=1,\ldots,m,\\
		~&~ 0\leq x_i\leq 1,~i=1,\ldots,n,
	\end{aligned}\right.
\end{equation}
where the feasible domain is assumed to be nonempty. Quadratic programming problems with linear inequality constraints can be equivalently expressed in this form by introducing slack variables.

Since (LCQP) is a special case of (QCQP), it admits a QNR as:
\begin{equation}\label{LCQPQNR0}
	\left\{\begin{aligned}
		\min ~&~ \frac{1}{2}x^{T}(Q-Z) x +c^Tx  +\frac{1}{2}t\\
		\mbox{s.t.} ~&~ a_i^T x = d_i,~i=1,\ldots,m,\\
		~&~ 0\leq x_i\leq 1,~i=1,\ldots,n,\\
		~&~ t= x^T Zx,
	\end{aligned}\right.
\end{equation}
where $Z\in\mathbb{S}^n$ satisfies $Q-Z\succeq 0$. The optimal parameter $Z^\ast$, which maximizes the McCormick relaxation bound of this reformulation, can be obtained via Proposition~\ref{thm1}.

However, this formulation does not fully leverage the linear equality constraints. A stronger reformulation arises by perturbing the objective function with the quadratic equalities $(a^T_i x-d_i)^2=0$ for $i=1,\ldots,m$, leading to:
\begin{equation}\label{LCQP-QNR}
	(\textrm{LCQP-QNR})~\left\{\begin{aligned}
		\min ~&~ \frac{1}{2}x^{T}Q x +c^Tx+\sum_{i=1}^m \gamma_i(a_i^T x-d_i)^2- \frac{1}{2}x^{T}Z x +\frac{1}{2}t\\
		\mbox{s.t.} ~&~ a_i^T x = d_i,~i=1,\ldots,m,\\
		~&~ 0\leq x_i\leq 1,~i=1,\ldots,n,\\
		~&~ t= x^T Zx,
	\end{aligned}\right.
\end{equation}
where $\gamma\in\mathbb{R}^m$ and $Z\in\mathbb{S}^n$ satisfy $Q+ 2\sum_{i=1}^m \gamma_i a_i a_i^T-Z\succeq 0$.
By incorporating the quadratic equalities, (LCQP-QNR) typically yields a tighter McCormick relaxation than the basic reformulation \eqref{LCQPQNR0}.

The approach of perturbing quadratic objective functions using $(a^T_ix-d_i)^2$ has been employed previously in quadratic convex reformulations \citep[Section 2]{Billionnet2012}. We extend this technique within the QNR framework. We use $\theta^\prime(\gamma,Z)$ to denote the McCormick relaxation bound of \eqref{LCQP-QNR}. The problem of finding the best parameters for constructing (LCQP-QNR) can be formulated as follows:
\begin{equation}\label{best_par2}
	\left\{\begin{aligned}
		\max ~&~ \theta^\prime(\gamma,Z) \\
		\mbox{s.t.} ~&~ Q+2\sum_{i=1}^m \gamma_i a_i a_i^T-Z\succeq 0,~\gamma\in\mathbb{R}^m.
	\end{aligned}\right.
\end{equation}
The following proposition formalizes our result:

\begin{proposition}\label{thm2}
	Assume that the optimal solution of \eqref{best_par2} is attainable.
	The optimal solution $(\gamma^\ast, Z^\ast)$ of \eqref{best_par2} can be obtained by solving the  semidefinite programming problem:
	\begin{equation}\label{QNR-SDP3}
		(\textrm{QNR-LCQP-SDP})~\left\{\begin{aligned}
			\max ~&~ \tau \\
			\textrm{s.t.} ~&~\begin{bmatrix}
				2\hat{b}-2\tau &\hat{c}^T\\
				\hat{c} ~&\hat{Q} \\
			\end{bmatrix}\succeq 0,\\
			~&~\hat{Q}= Q+2\sum_{i=1}^m \gamma_i a_i a_i^T +2\textrm{diag}(\beta)-2M-2N+2R+2S,\\
			~&~\hat{c}=c-2\sum_{i=1}^m \gamma_i d_i a_i +\sum_{i=1}^m \mu_i a_i-\beta+2Ne-Re-Se,\\
			~&~\hat{b}=\sum_{i=1}^m \gamma_i d_i^2-\sum_{i=1}^m \mu_i d_i-\sum_{i=1}^n\sum_{j=1}^n N_{ij},\\
			~&~\beta\in \mathbb{R}^n_+,~\gamma,\mu\in \mathbb{R}^m, ~\tau\in\mathbb{R}, ~M,N,R,S\in\mathbb{X}_+^n,
		\end{aligned}\right.
	\end{equation}
	and setting
	\begin{equation}\label{recoverZ2}
		Z^\ast= -2\textrm{diag}(\beta^\ast)+2M^\ast+2N^\ast-2R^\ast-2S^\ast,
	\end{equation}
	where $(\gamma^\ast,\beta^\ast,M^\ast,N^\ast,R^\ast,S^\ast)$ is the optimal solution of (QNR-LCQP-SDP).
\end{proposition}

The proof follows similar lines as Proposition~\ref{thm1} and arguments from \citep[Section 2]{Billionnet2012} and is omitted here for brevity.

The dual problem of (QNR-LCQP-SDP) is given as follows:

\begin{equation}\label{LINSDPRLT}
	(\textrm{LCQP-SDP+RLT})~\left\{\begin{aligned}
		\min ~&~ \frac{1}{2}Q\cdot X+c^{T}x \\
		\mbox{s.t.} ~&~ (a_ia_i^T)\cdot X-2d_i a_i^{T}x+d_i^2=0,~i=1,\ldots,m,\\
		&~a_i^T x=d_i,~i=1,\ldots,m,\\
		&~ X_{ii}- x_i\leq 0,~i=1,\ldots,n,\\
		&~ (X_{ij},x_i,x_j)\in \mathcal{M}_{ij},~i,j=1,\ldots,n,~i\neq j,\\
		&~ X\succeq xx^T.
	\end{aligned}\right.
\end{equation}


Notice that if  $(\gamma^\ast,Z^\ast)$ is an optimal solution of \eqref{best_par2}, so is $(\gamma^\prime,Z^\ast)$ for any $\gamma^\prime\in \mathbb{R}^m$ such that $\gamma_i^\prime\geq \gamma^\ast_i$ for $i=1,\ldots,m$. Such an observation is also discussed in Section 2 of  \citep{Billionnet2012}. Thus, under the assumption that \eqref{best_par2} has an attainable solution, we may fix the entries of $\gamma$ to a large enough constant. 

In addition, we remark that under the assumption that (LCQP) is feasible, both (QNR-LCQP-SDP) and (LCQP-SDP+RLT) are feasible. Furthermore, we can check that (QNR-LCQP-SDP) is always strictly feasible. Thus, the gap between (LCQP-SDP+RLT) and (QNR-LCQP-SDP) is zero. However, note that under the constraints $(a_ia_i^T)\cdot X-2d_i a_i^{T}x+d_i^2=0$, which is equivalent to $(a_ia_i^T)\cdot (X-xx^T) +(a_i^{T}x-d_i)^2=0$,
problem (LCQP-SDP+RLT) does not have a strictly feasible solution (a feasible solution that satisfies $X\succ xx^T$) under the assumption that $m\geq 1$ and $a_i\neq 0$ for $i=1,\ldots,m$. Thus, the existence of an optimal solution of (QNR-LCQP-SDP) is not theoretically guaranteed. Solving
(QNR-LCQP-SDP) may fail to obtain an optimal solution.

To derive a practical quadratic nonconvex reformulation with theoretical guarantee, we solve (LCQP-QNR) approximately rather than exactly by fixing each entry of $\gamma$ to a large enough constant. We use (QNR-LCQP-SDP$_{\gamma}$) to denote the reduced version of (QNR-LCQP-SDP) by regarding $\gamma$ as a constant vector, the dual problem of (QNR-LCQP-SDP$_{\gamma}$) is given as follows:

\begin{equation}\label{LINSDPRLT2}
	(\textrm{LCQP-SDP+RLT}_\gamma)~\left\{\begin{aligned}
		\min ~&~ \frac{1}{2}Q\cdot X+c^{T}x+\sum_{i=1}^m \gamma_i\left((a_ia_i^T)\cdot X-2d_i a_i^{T}x+d_i^2\right)\\
		\mbox{s.t.} ~&~a_i^T x=d_i,~i=1,\ldots,m,\\
		&~X_{ii}- x_i\leq 0,~i=1,\ldots,n,\\
		&~ (X_{ij},x_i,x_j)\in \mathcal{M}_{ij},~i,j=1,\ldots,n,~i\neq j,\\
		&~ X\succeq xx^T.
	\end{aligned}\right.
\end{equation}

We can see that both (QNR-LCQP-SDP$_\gamma$) and (LCQP-SDP+RLT$_\gamma$) are strictly feasible under some mild conditions that are generally satisfied in practice (e.g., there exists a feasible solution $x$ of (LCQP) that satisfies $0<x_i<1$ for $i=1,\ldots,n$). Thus, under these conditions, it is theoretically guaranteed that both the two problems have attainable optimal solutions and can be solved in polynomial-time. Fixing the entries of $\gamma$ to a large enough constant could derive a quadratic nonconvex reformulation with a high quality McCormick relaxation bound.

\textbf{Remark 3.} Although (LCQP-SDP+RLT) is not strictly feasible, (QNR-LCQP-SDP) may still have an attainable optimal solution. In our numerical experiments, we use Mosek, an interior-point algorithm-based solver,  to solve (QNR-LCQP-SDP) directly. If the solver fails to return an optimal solution, then we fix the entries of $\gamma$ to a large enough constant (we set it to be $10^4$ in our experiment), and solve (QNR-LCQP-SDP$_\gamma$) to obtain a sub-optimal solution of \eqref{best_par2}.

\section{Applications of QNR to classical quadratic optimization problems}\label{sec4}

In this section, we demonstrate the application of the proposed QNR technique to several well-known quadratic optimization problems.

\textbf{Example 1.} The first example is the box-constrained quadratic programming problem (BoxQP):
\begin{equation}\label{BoxQP}
	(\textrm{BoxQP})~\left\{\begin{aligned}
		\min ~&~\frac{1}{2}x^{T}Q x \\
		\mbox{s.t.} ~&~ 0\leq x_i \leq 1,~i=1,\ldots,n,
	\end{aligned}\right.
\end{equation}
where $Q\in\mathbb{S}^n$. This is a special case of (LCQP) without linear equality constraints. Following Section~\ref{sec3}, the QNR of (BoxQP) is:
\begin{equation}\label{BoxQP-QNR}
	(\textrm{BoxQP-QNR})~\left\{\begin{aligned}
		\min ~&~ \frac{1}{2}x^{T}Q x +c^Tx- \frac{1}{2}x^{T}Z x +\frac{1}{2}t\\
		\mbox{s.t.} ~&~ 0\leq x_i\leq 1,~i=1,\ldots,n,\\
		~&~ t= x^T Zx,
	\end{aligned}\right.
\end{equation}
where $Z\in \mathbb{S}^n$ satisfies $Q-Z\succeq 0$.

\textbf{Example 2.} The second example is the standard quadratic programming problem (StQP):
\begin{equation}\label{StQP}
	(\textrm{StQP})~\left\{\begin{aligned}
		\min ~&~\frac{1}{2}x^{T}Q x \\
		\mbox{s.t.} ~&~ \sum_{i=1}^n x_i=1,~x\geq 0.
	\end{aligned}\right.
\end{equation}
Since any feasible solution $x$ of (StQP) satisfies $0\leq x_i\leq 1$ for all $i=1,\ldots,n$, (StQP) is a special case of (LCQP). The QNR of (StQP) is:
\begin{equation}\label{StQP-QNR}
	(\textrm{StQP-QNR})~
	\left\{\begin{aligned}
		\min ~&~ \frac{1}{2}x^{T}(Q-Z) x+\gamma( \sum_{i=1}^n x_i-1)^2 +\frac{1}{2}t \\
		\mbox{s.t.} ~&~ \sum_{i=1}^n x_i=1,~t= x^T Zx,\\
		~&~ 0\leq x_i\leq 1,~i=1,\ldots,n.
	\end{aligned}\right.
\end{equation}
Note that although the upper bounds $x_i \le 1~(i=1,\ldots,n)$ are redundant, including them explicitly in the formulation is crucial for solvers to construct tight McCormick relaxations.

\textbf{Example 3.} The third example is the quadratically constrained quadratic program without explicit bounds:
\begin{equation}\label{QCQP2}
	(\textrm{QCQP2})~\left\{\begin{aligned}
		\min ~&~\frac{1}{2}x^{T}Q_0 x+c_0^{T}x \\
		\mbox{s.t.} ~&~\frac{1}{2}x^{T}Q_i x+c_i^{T}x\leq b_i,~i=1,2,\ldots,m,
	\end{aligned}\right.
\end{equation}
where the feasible domain is assumed bounded, but explicit variable bounds are not given.

To estimate variable bounds, we may solve the following bounding problems for each $x_i$:
\begin{equation}\label{Bounding}
	(\textrm{Bounding})~\left\{\begin{aligned}
		\min~(\textrm{or } \max): ~&~ x_i \\
		\mbox{s.t.} ~&~ \frac{1}{2}Q_i \cdot X+c_i^{T}x\leq b_i,~i=1,2,\ldots,m,\\
		~&~ X\succeq xx^T.
	\end{aligned}\right.
\end{equation}
Under the assumption that there exist $\gamma_1,\ldots,\gamma_m\geq 0$ such that $\sum_{i=1}^m \gamma_i Q_i \succ 0$, Slater's condition holds for \eqref{Bounding}, guaranteeing finite lower and upper bounds $l_i$ and $u_i$ can be obtained by solving  \eqref{Bounding}.

Adding these bounds, $l_i\leq x_i\leq u_i~(i=1,\ldots,n)$, transforms (QCQP2) into a special case of (QCQP), allowing us to apply the QNR. Like in the previous case, although these bounds are redundant, explicitly including them is essential to construct strong McCormick relaxations.

\section{Numerical Results}\label{sec5}

This section presents numerical results demonstrating the effectiveness of the proposed QNR technique. We conducted tests on four classes of problems: standard quadratic programming (StQP), box-constrained quadratic programming (BoxQP), general linearly constrained quadratic programming (LCQP), and general nonconvex quadratically constrained quadratic programming (QCQP).

To assess QNR's impact, we compare two approaches:
\begin{itemize}
	\item[$\bullet$]\textbf{Gurobi:} Solves the original problem formulations using Gurobi's default settings without any reformulation.
	\item[$\bullet$]\textbf{QNR:} Solves the QNR-reformulated problem using Gurobi, with presolve disabled  (by setting the parameter Presolve=0 in Gurobi) to preserve the QNR structure.
\end{itemize}

All experiments were run on a personal computer equipped with an AMD Ryzen 9 7940H CPU (4.00 GHz) and 16 GB RAM. The QNR related preprocessing were implemented in MATLAB R2023b. Semidefinite programming (SDP) problems were solved using Mosek 10.2.2, while the final reformulated problems were solved with Gurobi 11.02.

For fairness, all tests used a single computational thread and a relative optimality tolerance of $10^{-4}$. Other solver parameters remained at their default values.

\textbf{Remark 4.} Although Gurobi offers multiple lower-bounding techniques, we verified that for all tested instances, the root-node bounds computed by Gurobi coincided with the McCormick relaxation bounds. This confirms that McCormick relaxation was the active default bounding method. Details of this verification are provided in Appendix B.

\textbf{Remark 5.} Additional experiments using SCIP were conducted to verify the broader applicability of QNR. These results are reported in Appendix C.

Constructing the QNR requires solving the semidefinite program (QNR-SDP), whose dual problem (SDP+RLT) can be computationally demanding due to the large number of McCormick inequalities. To mitigate this, we apply the following iterative relaxation procedure:
\begin{enumerate}
	\item Solve (SDP+RLT) initially with all McCormick inequalities on off-diagonal entries $X_{ij}$  $(i\neq j)$ being dropped.
	\item Identify violated McCormick inequalities and add them back into the relaxation.
	\item Re-solve the relaxation and extract dual solutions to obtain a near-optimal solution to (QNR-SDP).
\end{enumerate}

This approach balances the trade-off between bound quality and preprocessing time and was applied to all problem classes except standard quadratic programming problem.

For standard quadratic programming problem, we simplify by directly solving a reduced version of (LCQP-SDP+RLT) that only retains the nonnegativity constraints $X_{ij}\geq 0$ in $(X_{ij},x_i,x_j)\in \mathcal{M}_{ij}$, foregoing iterative refinement. This strategy effectively balances QNR construction time against branch-and-bound efficiency.

\subsection{Standard quadratic programming problems}
We first evaluate the proposed QNR technique on the standard quadratic programming problem (StQP) introduced in Section~\ref{sec4}. Numerical experiments were conducted on two publicly available benchmark sets, SQP30 and SQP50 \citep{Bonami2016}, as well as on a set of challenging randomly generated instances created following the procedure in \citep{Gokmen}.

In addition to Gurobi and QNR, we tested quadprogIP, a state-of-the-art solver developed by \cite{Xia2020}. QuadprogIP reformulates linearly constrained nonconvex quadratic programs as mixed-integer linear programs and solves them with CPLEX. Prior studies by \cite{Xia2020} have shown that quadprogIP outperforms quadprogBB \citep{Chen} and other solvers, making it a strong benchmark for comparison.

Each of the SQP30 and SQP50 sets contains 150 instances with $n=30$ and $n=50$ variables, respectively. All instances were solved using the three methods with a 1200-second time limit.

\begin{figure}[htbp]
	\centering
	\begin{minipage}{1.0\linewidth}
		\centering
		\includegraphics[width=1.0\linewidth]{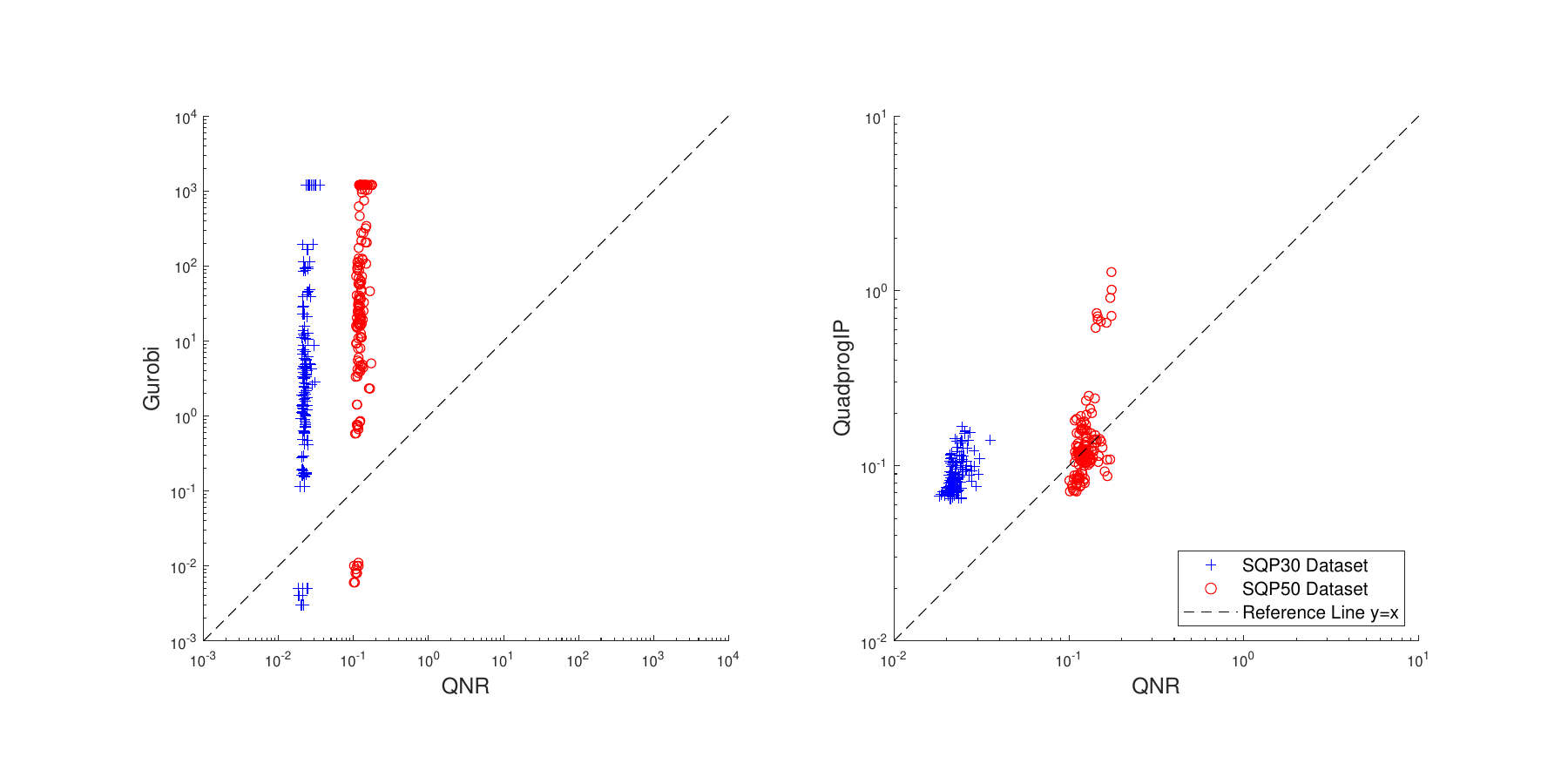}
	\end{minipage}
	\vspace{2mm}
	\caption{\centering{Log-log plots comparing QNR with Gurobi and quadprogIP on SQP30 and SQP50 test sets.}}\label{fig_stQP}
\end{figure}


Figure~\ref{fig_stQP} presents log-log plots comparing computational times of QNR against Gurobi and quadprogIP, where QNR's runtime includes both the time to solve the semidefinite program for reformulation and the subsequent branch-and-bound phase.

The left plot in Figure~\ref{fig_stQP} shows that QNR outperforms Gurobi on the majority of instances. The worst-case computation times for QNR were 0.04 seconds (SQP30) and 0.18 seconds (SQP50), with mean times of 0.02 and 0.13 seconds on the two sets, respectively. In contrast, Gurobi failed to solve 10 instances in SQP30 and 35 instances in SQP50 within the 1200-second limit, and its average times on solved instances were 16.85 and 89.47 seconds, respectively. These results highlight substantial runtime improvements enabled by QNR.

The right plot shows that QNR outperforms quadprogIP on all SQP30 instances, and achieves similar efficiency as quadprogIP on most SQP50 instances. In addition, we can see that there are ten instances in SQP50 set that are solved by quadprogIP using more than 0.5 seconds, while these instances are solved by QNR within 0.18 seconds. For quadprogIP, worst-case times were 0.17 seconds (SQP30) and 1.28 seconds (SQP50), with mean times of 0.09 and 0.17 seconds, respectively. These findings indicate superior performance of QNR on these benchmark sets.

We note that the SQP30 and SQP50 instances are relatively easy for both QNR and quadprogIP, as both solvers quickly solve all instances, limiting discriminatory power. Indeed, for all SQP30 and SQP50 instances, QNR achieves a zero relaxation gap at the root node, enabling solution without branching.

Other publicly available test sets used in \citep{Xia2020}--namely SQP, StableQP, and Scozzari/Tardella--were also tested. Most instances in these sets are solved by QNR in the root node. Therefore, to provide a more rigorous comparison, we generated \textbf{hard instances} following the procedure proposed by \cite{Gokmen}, summarized below in two steps.
\begin{itemize}
	\item[$\bullet$] \textbf{Step 1: Generate Matrix $\widehat{Q}$:} For a given $n\geq 5$, define the matrix $\widehat{Q} \in \mathbb{S}^n$ as $$\widehat{Q} \overset{\mathrm{def}}{=} \begin{bmatrix} B & C \\ C^T & H \end{bmatrix},$$ where
	\begin{equation*}
		H \overset{\mathrm{def}}{=} \begin{bmatrix}
			1 & -1 & 1 & 1 & -1 \\
			-1 & 1 & -1 & 1 & 1 \\
			1 & -1 & 1 & -1 & 1 \\
			1 & 1 & -1 & 1 & -1 \\
			-1 & 1 & 1 & -1 & 1
		\end{bmatrix}.
	\end{equation*}
	Matrix $B\in\mathbb{S}^{n-5}$ is a randomly generated positive semidefinite matrix of order $n-5$, and $C\in\mathbb{R}^{(n-5)\times 5}_+$ is a randomly generated matrix with nonnegative entries. In our experiments, the matrix $B$ is set to $B=VV^T$, where each element of $V\in\mathbb{R}^{(n-5)\times (n-5)}$ is uniformly sampled from interval $[-50,50]$, and each element of $C$ is uniformly sampled from $[0,1]$.
	\item[$\bullet$] \textbf{Step 2: Generate matrix $Q$:} Set $Q = J D \widehat{Q} D J^T$, where $J\in\mathbb{R}^{n\times n}$ is an arbitrary permutation matrix, and  $D\in\mathbb{S}^n$ is a diagonal matrix with entries sampled uniformly from $[0,1]$.
\end{itemize}

Using this procedure, we created 55 challenging instances for $n\in\{5,10,...,50,55\}$ (five instances per size), all exhibiting nonzero doubly nonnegative relaxation gaps \citep{Gokmen}.

We solved these hard instances using Gurobi, QNR, and quadprogIP with a 1200-second time limit. Table~\ref{tab1} summarizes the computational results, where each row lists the average number of solved instances (``Sol'') and the average runtime on the solved instances in seconds (``Time'') per method, grouped by problem size $n$.

\begin{table}[h]
	\caption{Results on randomly generated hard standard quadratic programming test instances}\label{tab1}
	\centering
	\begin{tabular}{lcccccc}
		\toprule
		\multirow{2}*{$n$} & \multicolumn{2}{c}{Gurobi} & \multicolumn{2}{c}{quadprogIP} & \multicolumn{2}{c}{QNR}  \\
		\cmidrule(lr){2-3} \cmidrule(lr){4-5} \cmidrule(lr){6-7}
		& Sol & Time & Sol & Time & Sol & Time  \\
		\midrule
		5  &5 &4.01    &5 &0.10  &5 &0.09  \\
		10  &5 &33.04    &5 &0.07  &5 &0.08  \\
		15  &2 &430.95    &5 &0.08  &5 &0.25  \\
		20  &0 &-    &5 &0.11  &5 &0.45  \\
		25  &0 &-    &5 &0.25  &5 &0.87  \\
		30  &0 &-    &5 &1.93  &5 &0.63  \\
		35  &0 &-    &5 &105.12  &5 &1.47  \\
		40  &0 &-    &4 &5.86  &5 &1.83  \\
		45  &0 &-    &5 &135.05  &5 &1.26  \\
		50  &0 &-    &4 &54.37  &5 &2.21  \\
		55  &0 &-    &4 &370.73  &5 &3.07  \\	
		\bottomrule
	\end{tabular}
\end{table}


Table 1 shows that: (1) For $n\leq 25$, quadprogIP outperformed QNR, but its performance degraded significantly as $n$ increased. (2) When $n\geq 35$, QNR was substantially faster, solving all instances within seconds, while quadprogIP failed to solve several instances for $n\geq 40$. (3) Gurobi, solving the original formulation, only managed instances with $n\leq 15$ and failed on all larger instances. In contrast, Gurobi combined with QNR solved all instances efficiently.

These results demonstrate that QNR markedly improves the performance of Gurobi on challenging instances, especially where Gurobi alone struggles.

\subsection{Box-constrained quadratic programming problems}

We further evaluate the QNR technique on the box-constrained quadratic programming problem (BoxQP). Our test instances come from three publicly available benchmark sets: the Basic Set \citep{Vandenbussche2005}, the Extended Set \citep{Burer2009}, and the Large Set \citep{Burer2010}, totaling 99 instances. Preliminary analysis revealed significant variation in the performance of Gurobi with respect to problem size and density. To better understand these effects, we partitioned the instances into two groups -- \textbf{Group A}: 81 instances with up to 90 variables, and \textbf{Group B}: 18 instances with 100 or more variables.

In addition to comparing QNR with Gurobi, we include FixingBB \citep{Locatelli2024} as a benchmark, since it is currently the fastest solver for globally solving BoxQP. We also evaluated QuadprogBB \citep{Chen}, a former state-of-the-art solver, but it was consistently slower than Gurobi and thus excluded from further comparisons. All three solvers--Gurobi, QNR, and FixingBB--were run on all 99 instances with a time limit of 7200 seconds, following the setup in \citep{Locatelli2024}.

\begin{figure}[b]
	\centering
	\begin{minipage}{1.0\linewidth}
		\centering
		\includegraphics[width=1.0\linewidth]{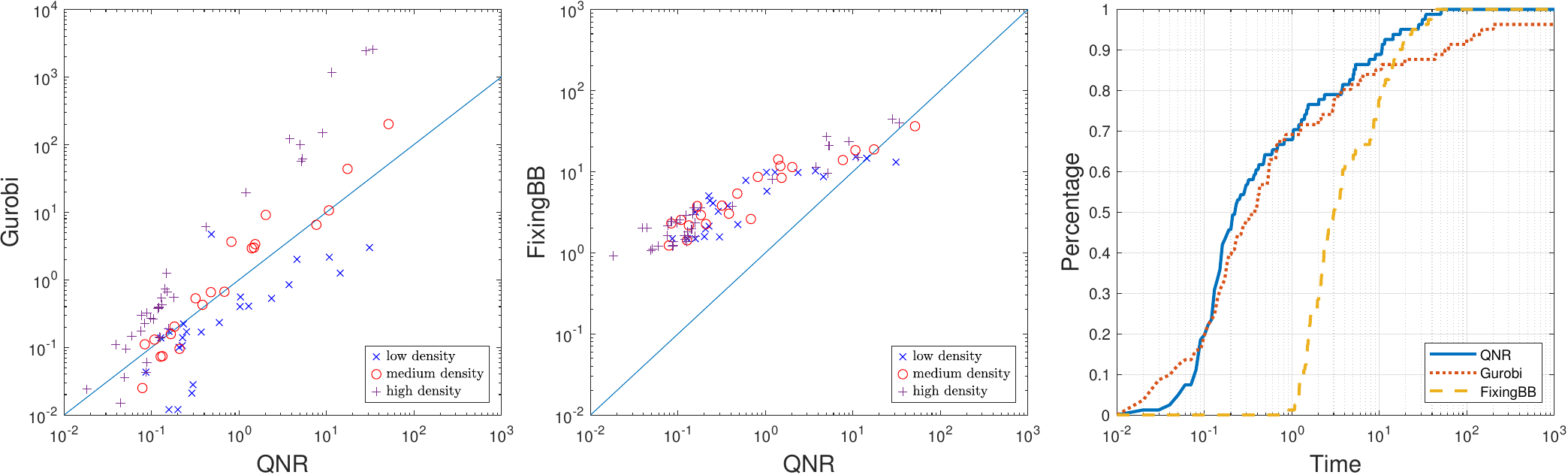}
	\end{minipage}
	\vspace{2mm}
	\caption{\centering{Log-log plots and performance profile for Group A}}\label{fig_GroupA}
\end{figure}

\begin{figure}[t]
	\centering
	\begin{minipage}{1.0\linewidth}
		\centering
		\includegraphics[width=1.0\linewidth]{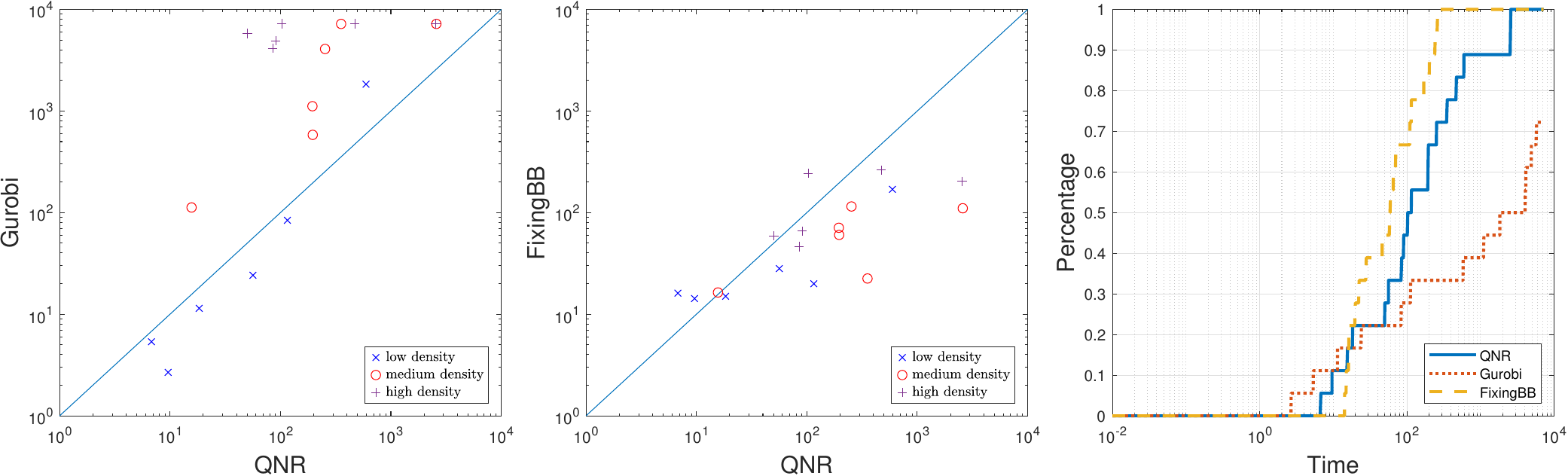}
	\end{minipage}
	\vspace{2mm}
	\caption{\centering{Log-log plots and performance profile for Group B}}\label{fig_GroupB}
\end{figure}

Figures~\ref{fig_GroupA} and \ref{fig_GroupB} present results via log-log plots and performance profiles, respectively, showing the percentage of problems solved within given time limits. To illustrate the impact of problem density, the log-log plots use different markers: $\times$ for low density ($d\leq 0.4$),  $\circ$ for medium density ($0.4<d\leq 0.6$), and $+$ for high density ($d>0.6$).

The two figures show that when QNR and Gurobi are compared, the effectiveness of QNR depends on problem density: (1) Among 30 low-density instances, QNR is slower than Gurobi on 26. (2) Among 27 medium-density instances, QNR outperforms Gurobi on 19 and under-performs on 8. (3) Among 42 high-density instances, QNR is faster on 38 and slower on only 4.

These findings indicate that QNR is particularly advantageous on medium- and high-density problems but less effective on sparse cases. This trend can be attributed to the McCormick relaxation: sparse problems require fewer McCormick inequalities since many $X_{ij}$ variables correspond to zero entries in $Q$, whereas QNR tends to introduce additional nonzero terms, increasing computational effort. Conversely, for denser problems, QNR's tighter McCormick relaxations significantly improve bound quality, resulting in better performance. Disregarding density, the performance profiles in Figures~\ref{fig_GroupA}~and~\ref{fig_GroupB} show that QNR consistently solves more instances than Gurobi when the runtime threshold exceeds 10 seconds, indicating that QNR is more effective on harder instances.

When QNR and FixingBB are compared, we see mixed results: (1) QNR outperforms FixingBB on 79 of the 81 instances in Group A. (2) FixingBB outperforms QNR on 13 of the 18 instances in Group B. Thus, we may conclude that QNR performs better on problems with up to 90 variables, while FixingBB dominates on larger instances with 100 or more variables.

This difference arises from their relaxation strategies: FixingBB applies an SDP+RLT relaxation at every node of the branch-and-bound tree, achieving very tight lower bounds at the expense of higher computational cost. By contrast, QNR applies the SDP+RLT relaxation only once at the root node, balancing relaxation strength and efficiency for moderate-sized problems. However, as the branch-and-bound tree grows, QNR's fixed parameters $Z$ do not update, causing lower bound tightness to degrade. For larger problems, where tight bounds throughout the enumeration are critical, FixingBB's approach yields better overall performance.

In summary, QNR delivers its greatest benefits on BoxQP instances with up to 90 variables and density above 0.4, significantly improving on Gurobi's performance and outperforming FixingBB in this regime. For larger problems, FixingBB's consistently tighter relaxations lead to superior performance. Additionally, QNR enhances Gurobi's efficiency on medium- and high-density problems but may increase computational complexity on very sparse instances due to the introduction of nonzero terms in the McCormick relaxation.

\subsection{General linearly constrained quadratic programming problems}

We next evaluate the proposed QNR technique on quadratic programming problems with linear inequality constraints, formulated as follows:
\begin{equation}\label{LCQP-INQ}
	(\textrm{LCQP-INQ})~\left\{\begin{aligned}
		\min ~&~\frac{1}{2}x^{T}Q x +c^Tx\\
		\mbox{s.t.} ~&~ a_i^T x\leq d_i,~i=1,\ldots,m,\\
		~&~ 0\leq x_i\leq 1,~i=1,\ldots,n.
	\end{aligned}\right.
\end{equation}
To apply QNR, we first convert (LCQP-INQ) into the standard (LCQP) form by introducing slack variables for the linear inequalities. This reformulated problem is then addressed using the QNR framework described in Section~\ref{sec3}.

We compare QNR's performance against two methods: Gurobi (without reformulation) and quadprogIP \citep{Xia2020}, which reformulates the problem as a mixed-integer linear program and solves it using CPlex. To simulate a variety of problem difficulties, we generate random instances for given $(m,n)$ pairs in the following manner: (1) Entries of $Q\in\mathbb{S}^n$ and $c\in\mathbb{R}^n$ are sampled uniformly from $[-10,10]$. (2) Each constraint vector $a_i\in \mathbb{R}^n$ has entries uniformly sampled from $[0,10]$. (3) Right-hand side values are set as $d_i = r_i\sum_{j=1}^n a_{ij}$, with $r_i$ sampled uniformly from $[0.2, 0.4]$. 

We start with the case of $m=1$. Test instances were generated with $n\in\{25,30,\ldots,60\}$, five per problem size. Table~\ref{result_qkp} summarizes results comparing Gurobi, quadprogIP, and QNR under a 1200-second time limit, reporting the number of instances solved and average runtime.

\begin{table}[h]
	\caption{Computational results on instances of (LCQP-INQ): the case of $m=1$}\label{result_qkp}
	\centering
	\begin{tabular}{crrrrrrrr}
		\toprule
		\multirow{2}*{$n$} & \multicolumn{2}{c}{Gurobi} &  \multicolumn{2}{c}{QuadprogIP} &  \multicolumn{2}{c}{QNR}  \\
		\cmidrule(lr){2-3} \cmidrule(lr){4-5} \cmidrule(lr){6-7}
		~ & Sol & Time & Sol & Time  & Sol & Time \\
		\midrule
		25 &5 &0.26  &5 &2.49  &5 &0.27 \\
		30 &5 &0.73  &5 &44.43  &5 &0.42 \\
		35 &5 &12.04  &4 &305.24  &5 &1.25 \\
		40 &5 &79.44  &0 &--  &5 &1.34 \\
		45 &5 &198.52  &0 &--  &5 &2.76 \\
		50 &3 &256.90  &0 &--  &5 &7.81 \\
		55 &0 &--  &0 &--  &5 &23.41 \\
		60 &1 &479.97  &0 &--  &5 &49.61 \\
		\bottomrule
	\end{tabular}
\end{table}

Table~\ref{result_qkp} shows that QNR consistently outperforms Gurobi for $n\geq 30$, solving more instances and with significantly shorter runtimes. The advantage grows with increasing $n$, demonstrating QNR's scalability and robustness. QNR also outperforms quadprogIP, which appears inefficient for this problem class. This inefficiency can be attributed to the reformulation it employs: quadprogIP reformulates the problem into a standard form defined in \citep{Xia2020} by introducing $m+n$ slack variables to convert the inequality constraints $a_i^T x\le d_i$ and $x_i \le 1$ into equality constraints, and then converts it into a mixed-integer linear programming problem. Such an expansion significantly increases the complexity of the reformulated problem, especially when $n$ is moderate to large. Thus, quadprogIP is not efficient for such problem, and runs slower than both Gurobi and QNR.

We then extend to experiments with varying $m$ and $n$. QNR's robustness are tested on instances with $m\in\{5,10,15\}$ and $n\in\{35,40,45,50\}$, generating five random instances per  $(n, m)$ pair. These instances are solved using Gurobi and QNR. QuadprogIP is not compared, since it is not efficient for such problem.  Table~\ref{result_lcqp} presents the results.

\begin{table}[h]
	\caption{Computational results on instances of (LCQP-INQ)}\label{result_lcqp}
	\centering
	\begin{tabular}{rrrrrrr}
		\toprule
		\multirow{2}*{$(n,m)$} & \multicolumn{2}{c}{Gurobi} &  \multicolumn{2}{c}{QNR}  \\
		\cmidrule(lr){2-3} \cmidrule(lr){4-5}
		~ & Sol & Time  & Sol & Time \\
		\midrule
		(35,5) &5 &36.61 &5 &1.34 \\
		(35,10) &5 &166.04 &5 &3.99 \\
		(35,15) &5 &93.39 &5 &1.98 \\
		(40,5) &5 &99.73 &5 &2.92 \\
		(40,10) &5 &163.33 &5 &1.77 \\
		(40,15) &3 &90.79 &5 &19.67 \\
		(45,5) &4 &182.74 &5 &24.86 \\
		(45,10) &1 &112.24 &5 &24.45 \\
		(45,15) &3 &586.88 &5 &77.08 \\
		(50,5) &1 &38.54 &5 &27.69 \\
		(50,10) &1 &159.11 &5 &106.38 \\
		(50,15) &0 &-- &5 &195.50 \\
		\bottomrule
	\end{tabular}
\end{table}

Table~\ref{result_lcqp} shows that QNR consistently achieves significantly faster runtimes than Gurobi across all configurations and solves all instances within the time limit. By contrast, Gurobi fails to solve several instances as $n$ and $m$ increase, particularly at $n=50$ and $m=15$. These results demonstrate that QNR effectively improves the efficiency of branch-and-bound algorithms in modern solvers like Gurobi, by tightening the McCormick relaxations.

\subsection{Quadratic programming problems with quadratic constraints}

We further assess the performance of the proposed QNR technique on a class of nonconvex quadratically constrained quadratic programming problems (QCQP). These problems are notoriously difficult to solve due to multiple nonconvex quadratic constraints, which substantially increase the complexity of relaxation and enumeration procedures in global optimization.

The testing instances are generated as following: Given a positive integer $m$ (the number of constraints), for each $i=0,1,\ldots,m-1$, we start by sampling a symmetric matrix $\hat{Q}_i$ with entries uniformly drawn from $[-10,10]$. We then perform eigenvalue decomposition $\hat{Q}_i= V_i^T \hat{D}_i V_i$, where $\hat{D}_i$ is diagonal and $V_i$ is orthogonal. Next, we construct matrix $Q_i = V_i^T D_i V_i$, where $D_0$ is diagonal with entries sampled uniformly from $[1, 20]$, while for $i=1,...,m-1$, the diagonal entries of $D_i$ are sampled from $[-10,10]$, introducing nonconvexity in the constraints.
The vector $c_0$ for the linear term in the objective function is set to the zero vector, and for each $i=1,\ldots,m-1$, the vectors $c_i$ has entries sampled uniformly from $[-10,10]$. Scalars $b_i$ are sampled from $[-10,-1]$. Additionally, the constraint $x^Tx\leq 1000$ is added as the $m$-th constraint to guarantee boundedness of the feasible domain. The feasiblity of all the generated instances have been checked.

Since the original formulation lacks explicit bounds on the variables--essential for constructing McCormick relaxations--we estimate valid lower and upper bounds for each variable by solving the bounding problem described in Section~\ref{sec4}. These bounds are incorporated into both the original and the reformulated problems prior to solution. 

We test three problem sizes: $n=m\in\{5,10,15\}$, corresponding to the largest sizes solvable within reasonable time by Gurobi in its current version. For each size, 10 random instances are generated and solved by Gurobi and QNR under a 1200-second time limit. Table~\ref{tab_qcqp_nonconvex} reports the number of instances solved within the limit and average runtime over solved instances.

\begin{table}[h]
	\caption{Computational results on instances of (QCQP).}\label{tab_qcqp_nonconvex}
	\centering
	\begin{tabular}{crrrrrr}
		\toprule
		$(n,m)$ & \multicolumn{2}{c}{Gorubi} &  \multicolumn{2}{c}{QNR}  \\
		\cmidrule(lr){2-3} \cmidrule(lr){4-5}
		& Sol & Time & Sol & Time\\
		\midrule
		(5,5) &10 &0.24 &10  &0.18 \\
		(10,10) &10 &74.69 &10  &39.26 \\
		(15,15) &0 &-- &1  &764.77 \\
		\bottomrule
	\end{tabular}
\end{table}

The results in Table~\ref{tab_qcqp_nonconvex} reveal the challenging nature of these nonconvex QCQPs even at small scales. Gurobi solves all instances for $n=5$ and $n=10$ but fails completely at $n=15$. In contrast, QNR matches Gurobi's performance for smaller sizes and manages to solve one instance at $n=15$. This marginal improvement supports the potential of QNR to accelerate the solver, even in the presence of high nonconvexity.

\section{Conclusions}\label{sec6}

This paper introduced a general framework for quadratic nonconvex reformulation (QNR), an unconventional yet very effective approach to enhancing the performance of branch-and-bound-based global solvers for nonconvex quadratic programming problems. By systematically tightening McCormick relaxations, QNR delivers substantial computational improvements for state-of-the-art solvers such as Gurobi and SCIP.

Extensive experiments across diverse benchmark sets confirm the versatility and effectiveness of QNR. For standard quadratic programming problems, QNR outperforms quadprogIP on challenging instances. For box-constrained quadratic programs, QNR outperforms the leading method, FixingBB, on public test sets with up to 90 variables. In general linearly constrained quadratic programs, integrating QNR into Gurobi yields marked performance gains.

The benefits are less pronounced for problems with nonconvex quadratic constraints, as QNR primarily reduces the relaxation gap of McCormick relaxations without directly tightening convex relaxations of the feasible region. Addressing this limitation presents a promising avenue for future research, particularly in extending QNR to exploit alternative relaxations or hybrid reformulation-relaxation strategies.

\section*{Appendix}\label{sec7}

This Appendix include the following parts: Appendix A presents an independent and detailed proof of Proposition 1 stated in Section~\ref{sec2}. This result provides theoretical justification for the proposed QNR technique. Appendix B offers additional numerical experiments that investigate how the proposed QNR technique influences the McCormick relaxation bounds computed by Gurobi. These results provide deeper insight into the impact of the reformulation on the performance of a solver. Appendix C reports supplementary computational results on the performance of SCIP when the QNR technique is applied. These results demonstrate the broad applicability of the proposed approach.


\subsection*{Appendix A: Proof for Proposition 1}\label{app1}
Although Proposition 1 in Section~\ref{sec2} can be viewed as an extension of earlier theoretical results--namely, \citet[Theorem 1]{Billionnet2016} and \citet[Theorem 1]{Zheng2011b}--certain technical details are different. For the sake of completeness and clarity, we provide an independent and self-contained proof of the proposition.

\begin{proof}{Proof}
	Recall the definition of the McCormick relaxation bound $\theta(Z_0, \ldots, Z_m)$ for the QNR reformulation, which corresponds to the optimal value of \eqref{QNR-MCR}.
	Define the Lagrangian function of \eqref{QNR-MCR} as
	\begin{equation}
		\begin{aligned}
			&L(x,X,\alpha,\beta,\gamma,M,N,R,S)=\frac{1}{2}x^{T}(Q_0-Z_0)x+c_0^{T}x+\frac{1}{2}Z_0\cdot X +\sum_{i=1}^n\alpha_i(x_i^2-X_{ii})\\
			&\quad +\sum_{i=1}^n\beta_i(X_{ii}-x_i)+\sum_{i=1}^m \gamma_i\left(\frac{1}{2}x^{T}(Q_i-Z_i)x+c_i^{T}x+ \frac{1}{2} Z_i\cdot X- b_i\right)\\
			&\quad +\sum_{i=1}^n\sum_{j=1}^n\left[ -M_{ij}X_{ij}-N_{ij}\left(X_{ij}-x_j-x_i+1\right) ) \right]\\
			&\quad +\sum_{i=1}^n\sum_{j=1}^n\left[  R_{ij}(X_{ij}-x_i)+S_{ij}(X_{ij}-x_j) \right],
		\end{aligned}
	\end{equation}
	where $\alpha,\beta\in \mathbb{R}^n_+$ and $\gamma\in \mathbb{R}^m_+$ are the Lagrangian multipliers of the corresponding inequality constraints in \eqref{QNR-MCR}, and $M, N, R, S \in \mathbb{X}^n_+$ whose entries $M_{ij}, N_{ij}, R_{ij}, S_{ij}$ are the Lagrangian multipliers of the four McCormick inequalities in the constraint $\left(X_{ij},x_i,x_j\right)\in \mathcal{M}_{ij}$, respectively. Define the dual function:
	\begin{equation}\label{dualfunction}
		p(\alpha,\beta,\gamma,M,N,R,S)=\min_{x\in\mathbb{R}^n,X\in\mathbb{S}^n}L(x,X,\alpha,\beta,\gamma,M,N,R,S),
	\end{equation}
	and the associated dual problem of \eqref{QNR-MCR}:
	\begin{equation}\label{QNR-dual}
		(\textrm{QNR-dual})~\left\{\begin{aligned}
			\max ~&p(\alpha,\beta,\gamma,M,N,R,S)\\
			\mbox{s.t.} ~~&\alpha,\beta\in \mathbb{R}^n_+,~\gamma\in \mathbb{R}^m_+, ~M,N,R,S\in\mathbb{X}_+^n.
		\end{aligned}\right.
	\end{equation}
	
	To avoid the unboundedness in minimizing the Lagrangian function in \eqref{dualfunction}, the coefficient for the matrix $X$ in $L(x,X,\alpha,\beta,\gamma,M,N,R,S)$ must vanish. This leads to the following necessary constraint:
	\begin{equation}\label{eq_condition}
		\frac{1}{2}Z_0+\sum_{i=1}^m\frac{1}{2} \gamma_i Z_i-\textrm{diag}(\alpha)+\textrm{diag}(\beta)-M-N+R+S= 0.
	\end{equation}
	Under this condition, the Lagrangian function simplifies to
	\begin{equation}
		\begin{aligned}
			L(x) = \frac{1}{2}x^T \hat{Q} x + \hat{c}^T x + \hat{b},
		\end{aligned}
	\end{equation}
	where
	\begin{equation}\label{Qcb}
		\begin{aligned}
			&\hat{Q}=(Q_0-Z_0)+\sum_{i=1}^m \gamma_i (Q_i-Z_i) +2\textrm{diag}(\alpha),\\
			&\hat{c}=c_0+\sum_{i=1}^m \gamma_i c_i -\beta+2Ne-Re-Se,\\
			&\hat{b}=-\sum_{i=1}^m \gamma_i b_i-\sum_{i=1}^n\sum_{j=1}^n N_{ij}.
		\end{aligned}
	\end{equation}
	Since $\gamma\geq 0$, $\alpha\geq 0$, and $Q_i-Z_i\succeq 0$ for all $i=0,\ldots,m$, we have $\hat{Q}\succeq 0$. By substituting the constraint \eqref{eq_condition}, the expression for $\hat{Q}$ can equivalently be written as
	\begin{equation}
		\hat{Q}=Q_0+\sum_{i=1}^m\gamma_i Q_i +2\textrm{diag}(\beta)-2M-2N+2R+2S.
	\end{equation}
	The dual problem can now be recast as the following program:
	\begin{equation}\label{QNR-dual2}
		\left\{\begin{aligned}
			\max ~&~ \tau \\
			\mbox{s.t.} ~&~L(x,X,\alpha,\beta,\gamma,M,N,R,S)-\tau\geq 0,~\forall x\in\mathbb{R}^n,~\forall X\in\mathbb{X}^n,\\
			~&~\frac{1}{2}Z_0+\sum_{i=1}^m\frac{1}{2} \gamma_i Z_i-\textrm{diag}(\alpha)+\textrm{diag}(\beta)-M-N+R+S= 0,\\
			~&~ \alpha\in \mathbb{R}^n_+,~\beta\in \mathbb{R}^n_+,~\gamma\in \mathbb{R}^m_+, ~\tau\in\mathbb{R}, ~M,N,R,S\in\mathbb{X}_+^n.
		\end{aligned}\right.
	\end{equation}
	Besides, the constraint $L(x,X,\alpha,\beta,\gamma,M,N,R,S)-\tau\geq 0~(\forall x\in\mathbb{R}^n,~\forall X\in\mathbb{X}^n)$ is equivalent to
	\begin{equation}\label{sdpconstraint}
		\begin{bmatrix}
			2\hat{b}-2\tau &\hat{c}^T\\
			\hat{c} ~&\hat{Q} \\
		\end{bmatrix}\succeq 0
	\end{equation}
	under the condition that $\hat{Q}\succeq 0$. Thus, we can replace the first constraint in \eqref{QNR-dual2} by \eqref{sdpconstraint}, and reformulate \eqref{QNR-dual2} as a semidefinite programming problem.
	
	By strong duality, which holds under the assumption of strict feasibility of (QCQP), we have $\theta(Z_0, \ldots, Z_m)$ equal to the optimal value of \eqref{QNR-dual2}. Thus, problem \eqref{best_par} can be reformulated by using \eqref{QNR-dual2} to represent $\theta(Z_0,\ldots,Z_m)$, which results in the following problem:
	\begin{equation}\label{QNR-SDP2}
		\left\{\begin{aligned}
			\max ~&~ \tau \\
			\mbox{s.t.} ~&~ \begin{bmatrix}
				2\hat{b}-2\tau &\hat{c}^T\\
				\hat{b} ~&\hat{Q} \\
			\end{bmatrix}\succeq 0,\\
			~&~ Q_i-Z_i\succeq 0,~i=0,\ldots,m,\\
			~&~\frac{1}{2}Z_0+\sum_{i=1}^m\frac{1}{2} \gamma_i Z_i-\textrm{diag}(\alpha)+\textrm{diag}(\beta)-M-N+R+S= 0,\\
			~&~\hat{Q}=Q_0+\sum_{i=1}^m\gamma_i Q_i +2\textrm{diag}(\beta)-2M-2N+2R+2S,\\
			~&~\hat{c}=c_0+\sum_{i=1}^m \gamma_i c_i -\beta+2Ne-Re-Se,\\
			~&~\hat{b}=-\sum_{i=1}^m \gamma_i b_i-\sum_{i=1}^n\sum_{j=1}^n N_{ij},\\
			~&~\alpha\in \mathbb{R}^n_+,~\beta\in \mathbb{R}^n_+,~\gamma\in \mathbb{R}^m_+, ~\tau\in\mathbb{R}, ~M,N,R,S\in\mathbb{X}_+^n,\\
			~&~Z_i\in\mathbb{S}^n,~i=0,\ldots,m.
		\end{aligned}\right.
	\end{equation}
	
	We now show that problem \eqref{QNR-SDP2} is equivalent to problem (QNR-SDP) defined in Section \ref{sec2}. Observe that problem (QNR-SDP) can be derived from \eqref{QNR-SDP2} by dropping the constraints $\frac{1}{2}Z_0+\sum_{i=1}^m\frac{1}{2} \gamma_i Z_i-\textrm{diag}(\alpha)+\textrm{diag}(\beta)-M-N+R+S= 0$ and $Q_i-2Z_i\succeq 0$ for $i=0,\ldots,m$ directly. Thus, the optimal value of (QNR-SDP) is no more than that of \eqref{QNR-SDP2}. On the other hand, for any feasible solution $(\gamma, \beta, M, N, R, S)$ of problem (QNR-SDP), we can recover a feasible solution of problem \eqref{QNR-SDP2} by setting
	\begin{equation}\label{recover}
		\begin{aligned}
			&\alpha=0,\\
			&Z_i=Q_i,~i=1,\ldots,m,\\
			&Z_0=-\sum_{i=1}^m \gamma_i Q_i-2\textrm{diag}(\beta)+2M+2N-2R-2S.
		\end{aligned}
	\end{equation}
	Thus, the optimal value of \eqref{QNR-SDP2} is no more than that of (QNR-SDP). Therefore, problem \eqref{best_par} and (QNR-SDP) are equivalent in the sense that they have the same optimal value, and the optimal solution of problem (6) can be obtained from the one of (QNR-SDP) via equation \eqref{LCQP-QNR}.
\end{proof}

\subsection*{Appendix B: Numerical Results on the McCormick Relaxation Bound in the Root Node} \label{app2}
In this section, we present computational results aimed at investigating the impact of the proposed QNR technique on the McCormick relaxation bounds utilized by Gurobi in the root node of the branch-and-bound framework.

Our experiments focus on two representative problem classes: (i) the linearly constrained quadratic programming problem with inequality constraints (LCQP-INQ), and (ii) the nonconvex quadratically constrained quadratic programming problem (QCQP), as described in Sections 5.3 and 5.4, respectively. For (LCQP-INQ), we select 15 instances with $n=40$ and $m\in\{5,10,15\}$, generated in Section 5.3. For (QCQP), we select 10 instances where $n = m = 10$, generated in Section 5.4. For each instance, we record the root node bounds returned by Gurobi, with and without the QNR reformulation, respectively. When applying the QNR technique, we solve the associated semidefinite relaxation problem (QNR-SDP) approximately using the iterative algorithm detailed in Section 5. These semidefinite relaxation bounds are also reported for comparison.

The results for the (LCQP-INQ) and (QCQP) instances are summarized in Tables~\ref{result_lcqp_app} and~\ref{result_qcqp_app}, respectively. The columns include the following:
\begin{itemize}
	\item[TG:] The time (in seconds) taken by Gurobi to solve the problem in the original or QNR formulation.
	\item[RB:] The root-node bound obtained by Gurobi in the original or QNR formulation.
	\item[TT:] Total time for QNR (the sum of the time for reformulation and running Gurobi).
	\item[SDR:] Semidefinite relaxation bound computed via solving problem (QNR-SDP).
	\item[OPT:] Global optimal value.
	\item[GC:] The gap closed by QNR, computed as
	\begin{equation}
		\textrm{GC}=1-\frac{\textrm{OPT}-\textrm{RB}_{\textrm{QNR}}}{ \textrm{OPT}-\textrm{RB}_{\textrm{Gurobi}} },
	\end{equation}
	where $\text{RB}_{\text{QNR}}$ and $\text{RB}_{\text{Gurobi}}$ denote the root node bounds with and without the QNR technique, respectively.
\end{itemize}

Each row corresponds to an individual instance, indexed by $k$. The symbol ``---'' indicates that the instance was not solved within the time limit of 1200 seconds.

Based on the numerical results reported in Tables~\ref{result_lcqp_app} and~\ref{result_qcqp_app}, we summarize the following observations:

\begin{enumerate}
	\item \textbf{Equivalence of QNR root node bound and SDR bounds.}
	In all test cases, the root node bound resulting from the QNR reformulation exactly matches the semidefinite relaxation bound. This result supports Proposition 1 and the theoretical discussion in Section 2.3, which states that problem \eqref{best_par} is equivalent to (QNR-SDP) in the sense that they have the same optimal value.
	
	\item \textbf{Significant tightening of relaxation bounds.}
	The QNR reformulation significantly tightens the root node bounds. For the 15 instances of problem (LCQP-INQ), over 94\% of the duality gap is closed. For the 10 instances of problem (QCQP), QNR closes more than 60\% of the gap on average. The tightening of the relaxation bounds explains why the proposed QNR technique can speed-up the solver based on the branch-and-bound framework.
	
	\item \textbf{Low overhead of QNR reformulation.}
	The time required to construct the QNR reformulation (which is equal to TT$-$TG) is small compared to the original Gurobi solution time. Consequently, the total solution time (including QNR construction) is often significantly less than that of solving the original problem directly.
	
	\item \textbf{Confirmation of McCormick relaxation as Gurobi's default relaxation.}
	Although the results are not listed in Tables~\ref{result_lcqp_app} and~\ref{result_qcqp_app}, we have verified that, across all instances considered in this paper, the root-node bound returned by Gurobi coincides with the McCormick relaxation bound (for both the cases of Gurobi and QNR). This confirms that Gurobi 11.02 defaults to using McCormick relaxations for handling nonconvexity in these instances.
\end{enumerate}

We remark that in Table~\ref{result_qcqp_app}, the root-node bounds returned by Gurobi are all zero. This is expected: for the test instances, the objective function is a homogeneous convex quadratic function, and  the McCormick relaxation of the problem always includes the origin $x=0$ in its feasible region. As a result, the default Gurobi relaxation returns a trivial lower bound of zero. In contrast, using the QNR technique, the structure of the problem is changed, and the root-node bound of McCormick relaxation can be strictly positive.

\begin{table}[h]
	\caption{Computational results on instances of (LCQP).}\label{result_lcqp_app}
	\centering
	\begin{tabular}{rrrrrrrrr}
		\toprule
		\multirow{2}*{$(n,m,k)$} & \multicolumn{2}{c}{Gurobi} &  \multicolumn{3}{c}{QNR} &\multirow{2}*{SDR} &\multirow{2}*{OPT} &\multirow{2}*{GC} \\
		\cmidrule(lr){2-3} \cmidrule(lr){4-6}
		~ &TG &RB  &TT &TG &RB &~  &~ &~ \\
		\midrule
		(40,5,1)  &1.56 &-402.94 &0.84 &0.53 &-209.20 &-209.20 &-208.90 &99.8\% \\
		(40,5,2)  &234.12 &-494.29 &4.55 &4.23 &-268.52 &-268.52 &-257.68 &95.4\% \\
		(40,5,3)  &113.61 &-406.22 &2.59 &2.27 &-169.69 &-169.69 &-161.63 &96.7\% \\
		(40,5,4)  &82.09 &-397.97 &4.20 &3.84 &-195.98 &-195.98 &-183.28 &94.1\% \\
		(40,5,5)  &67.30 &-388.78 &2.42 &2.09 &-194.06 &-194.06 &-182.83 &94.5\% \\
		(40,10,1)  &179.13 &-528.19 &2.28 &1.84 &-254.41 &-254.42 &-249.25 &98.1\% \\
		(40,10,2)  &507.71 &-433.33 &2.20 &1.81 &-184.34 &-184.34 &-175.76 &96.7\% \\
		(40,10,3)  &22.44 &-403.33 &1.30 &0.87 &-178.66 &-178.66 &-174.23 &98.1\% \\
		(40,10,4)  &46.65 &-424.74 &1.64 &1.23 &-215.90 &-215.90 &-209.82 &97.2\% \\
		(40,10,5)  &60.70 &-417.22 &1.45 &1.12 &-204.02 &-204.02 &-197.42 &97.0\% \\
		(40,15,1)  &261.73 &-510.31 &75.00 &74.40 &-251.33 &-251.33 &-237.41 &94.9\% \\
		(40,15,2)  &--- &-409.38 &3.08 &2.71 &-181.28 &-181.28 &-166.86 &94.1\% \\
		(40,15,3)  &--- &-435.27 &15.32 &14.88 &-187.32 &-187.33 &-171.98 &94.2\% \\
		(40,15,4)  &3.03 &-370.61 &3.00 &2.26 &-198.72 &-198.72 &-193.42 &97.0\% \\
		(40,15,5)  &7.60 &-418.06 &1.94 &1.44 &-218.22 &-218.22 &-215.22 &98.5\% \\
		\bottomrule
	\end{tabular}
\end{table}

\begin{table}[h]
	\caption{Computational results on instances of (QCQP).}\label{result_qcqp_app}
	\centering
	\begin{tabular}{rrrrrrrrr}
		\toprule
		\multirow{2}*{$(n,m,k)$} & \multicolumn{2}{c}{Gurobi} &  \multicolumn{3}{c}{QNR} &\multirow{2}*{SDR} &\multirow{2}*{OPT} &\multirow{2}*{GC} \\
		\cmidrule(lr){2-3} \cmidrule(lr){4-6}
		~ &TG &RB  &TT &TG &RB &~  &~ &~ \\
		\midrule
		(10,10,1)  &93.45 &0.00 &61.50 &61.50 &22.13 &22.13 &29.48 &75.1\% \\
		(10,10,2)  &59.16 &0.00 &0.96 &0.95 &12.53 &12.53 &12.55 &99.9\% \\
		(10,10,3)  &37.24 &0.00 &22.87 &22.87 &6.96 &6.96 &7.77 &89.6\% \\
		(10,10,4)  &89.53 &0.00 &25.18 &25.18 &7.91 &7.91 &8.81 &89.7\% \\
		(10,10,5)  &76.26 &0.00 &55.75 &55.74 &13.90 &13.90 &20.15 &69.0\% \\
		(10,10,6)  &104.07 &0.00 &74.28 &74.28 &12.34 &12.34 &20.40 &60.5\% \\
		(10,10,7)  &90.96 &0.00 &62.04 &62.04 &4.98 &4.98 &7.77 &64.2\% \\
		(10,10,8)  &36.27 &0.00 &1.91 &1.91 &5.27 &5.27 &5.30 &99.5\% \\
		(10,10,9)  &81.36 &0.00 &41.97 &41.97 &21.33 &21.33 &30.31 &70.4\% \\
		(10,10,10)  &78.59 &0.00 &46.14 &46.13 &13.43 &13.43 &17.92 &74.9\% \\
		\bottomrule
	\end{tabular}
\end{table}

\subsection*{Appendix C: Numerical Results of the QNR Technique on SCIP} \label{app3}

The applications of the proposed QNR technique are not limited to accelerating Gurobi. It can be potentially applied to any nonconvex quadratic programming solvers that employ the McCormick relaxation framework, as described in Section~2.1. In this section, we further investigate the impact of the QNR technique on SCIP, a well-known open-source solver for general mixed-integer nonlinear programming problems.

\subsubsection*{C1: Standard quadratic programming}

We first consider the standard quadratic programming problem using the same test instances described in Section~5.1. Note that for all instances in the benchmark sets SQP30 and SQP50, the QNR technique can yield reformulations with zero McCormick relaxation gaps. These reformulated instances can therefore be solved directly at the root node by both Gurobi and SCIP. To further evaluate the performance of QNR, we use the hard test instances randomly generated in Section 5.1.

Table~\ref{tab1_app} reports numerical results for SCIP and QNR+SCIP. Each row summarizes results for instances with a fixed dimension $n$, including the number of solved instances (column ``Sol'') within the 1200-second time limit and the average time (column ``Time'') required for solving them. The column ``SCIP'' presents the results using SCIP on the original formulation, while ``QNR+SCIP'' shows the results using SCIP on the QNR-reformulated instances.

\begin{table}[h]
	\caption{Results on randomly generated standard quadratic programming instances}\label{tab1_app}
	\centering
	\begin{tabular}{lcccc}
		\toprule
		\multirow{2}*{$n$} & \multicolumn{2}{c}{SCIP} & \multicolumn{2}{c}{QNR+SCIP}  \\
		\cmidrule(lr){2-3} \cmidrule(lr){4-5}
		& Sol & Time & Sol & Time   \\
		\midrule
		5 &5 &2.00  &5 &0.80 \\
		10 &0 &---  &5 &1.40 \\
		15 &0 &---  &5 &6.20 \\
		20 &0 &---  &5 &25.20 \\
		25 &0 &---  &5 &82.60 \\
		30 &0 &---  &5 &44.60 \\
		35 &0 &---  &5 &220.20 \\
		40 &0 &---  &4 &217.00 \\
		45 &0 &---  &5 &151.40 \\
		50 &0 &---  &2 &83.00 \\
		55 &0 &---  &1 &195.00 \\
		\bottomrule
	\end{tabular}
\end{table}

The results in Table~\ref{tab1_app} demonstrate that the QNR technique significantly enhances the performance of SCIP. Without QNR, SCIP fails to solve any instance with $n \geq 10$ within the time limit. By contrast, with QNR, SCIP successfully solves much larger test instances efficiently.

\subsubsection*{C2: Linearly inequality constrained quadratic programming}
Next, we evaluate the performance of the QNR technique on linearly inequality constrained quadratic programming problems (LCQPs). We use the  test instances generated in Section~5.3. The results are listed in Table~\ref{tab2_app}, where each row corresponds to instances of fixed dimensions $(n,m)$. The columns follow the same meaning as in Table~\ref{tab1_app}.
\begin{table}[h]
	\caption{Computational results on instances of (LCQP).}\label{tab2_app}
	\centering
	\begin{tabular}{rrrrrrr}
		\toprule
		\multirow{2}*{$(n,m)$} & \multicolumn{2}{c}{SCIP} &  \multicolumn{2}{c}{QNR+SCIP}  \\
		\cmidrule(lr){2-3} \cmidrule(lr){4-5}
		~ & Sol & Time  & Sol & Time \\
		\midrule
		(35,5) &3 &357.67  &5 &90.00 \\
		(35,10) &3 &440.00  &5 &403.20 \\
		(35,15) &5 &627.40  &5 &167.60 \\
		(40,5) &1 &560.00  &5 &228.60 \\
		(40,10) &0 &---  &5 &138.80 \\
		(40,15) &2 &413.50  &3 &277.67 \\
		(45,5) &0 &---  &3 &429.67 \\
		(45,10) &0 &---  &2 &685.00 \\
		(45,15) &0 &---  &3 &690.33 \\
		(50,5) &0 &---  &3 &675.00 \\
		(50,10) &0 &---  &3 &467.00 \\
		(50,15) &0 &---  &2 &365.50 \\
		\bottomrule
	\end{tabular}
\end{table}

From Table~\ref{tab2_app}, we observe that the QNR technique significantly improves the performance of SCIP across nearly all cases. SCIP without QNR fails to solve most of the instances with variable size $n\ge 40$ within 1200 seconds, whereas QNR+SCIP successfully solves many of them.

\subsubsection*{C3: Nonconvex quadratically constrained quadratic programming}
Finally, we evaluate the QNR technique on nonconvex quadratically constrained quadratic programming (QCQP) problems, using the test instances from Section~5.4. The results are reported in Table~\ref{tab_qcqp_nonconvex_app}.

From Table~\ref{tab_qcqp_nonconvex_app}, we find that the QNR technique provides limited acceleration when applied to SCIP on QCQP instances. For the smallest case $(5,5)$, QNR slightly increases the solution time. For $(10,10)$, the total solve time improves moderately. For $(15,15)$, both SCIP and QNR+SCIP fail to solve the instances within the time limit. It seems that the QNR technique is not effective in accelerating SCIP on nonconvex QCQP problems.

\begin{table}[h]
	\caption{Computational results on instances of (QCQP).}\label{tab_qcqp_nonconvex_app}
	\centering
	\begin{tabular}{crrrrrr}
		\toprule
		$(n,m)$ & \multicolumn{2}{c}{SCIP} &  \multicolumn{2}{c}{QNR+SCIP}  \\
		\cmidrule(lr){2-3} \cmidrule(lr){4-5}
		& Sol & Time & Sol & Time\\
		\midrule
		(5,5) &10 &2.7 &10  &3.9 \\
		(10,10) &7 &670.4 &7  &434.0 \\
		(15,15) &0 &-- &0  &-- \\
		\bottomrule
	\end{tabular}
\end{table}

The varying effectiveness of the QNR technique between (LCQP) and (QCQP) can be explained by the nature of McCormick relaxation. For (LCQP), the feasible region itself is convex and is not expanded by the McCormick relaxation. In such cases, the QNR technique effectively tightens the relaxation without altering the feasible domain (in the sense that the projection of the feasible domain to the space of $x$ keeps unchanged). In contrast, for nonconvex (QCQP), the McCormick relaxation relaxes the feasible domain, which introduces a more significant challenge. Although the QNR technique improves the relaxation of the objective function, it does not address the loose relaxation of the nonconvex constraints. Finding a feasible solution that satisfies all nonconvex quadratic constraints is already a NP-hard problem in general cases, and tightening the objective function alone is not sufficient to improve the performance of a solver.

Based on the above results, we conclude that the QNR technique is highly effective in accelerating SCIP on standard and linearly constrained quadratic programming problems. However, its benefits are limited when applied to nonconvex quadratically constrained quadratic programming problems.

\bibliographystyle{abbrvnat}
\bibliography{QNR}  






\end{document}